\documentclass[a4paper,12pt]{article}
\usepackage{stmaryrd}
\usepackage{bbm}
\usepackage{amsfonts}
\usepackage{amsmath}
\usepackage{amssymb}
\usepackage{mathrsfs}
\usepackage{accents}
\usepackage{amscd}

\numberwithin{equation}{section}

\title{\textbf{Twistor geometry of Riemannian 4-manifolds by moving frames}}
\author{Jixiang Fu~~~~Xianchao Zhou \footnote{E-mail addresses: majxfu@fudan.edu.cn,~~zhouxianch09@gmail.com}\\\\
Fudan University \\
School of Mathematical Sciences}
\date{}

\begin{document}

\maketitle

\noindent\textbf{Abstract.} In this paper, we characterize  Riemannian 4-manifold in terms of its
almost Hermitian twistor spaces $(Z,g_t,\mathbb{J}_{\pm})$. Some special metric conditions (including Balanced metric condition, first Gauduchon metric condition) on $(Z,g_t,\mathbb{J}_{\pm})$ are studied. For the first Chern form of a natural unitary connection on the vertical tangent bundle over the twistor space $Z$, we can recover J. Fine and D. Panov's result on the condition of the first Chern form being symplectic and P. Gauduchon's result on the condition of the first Chern form being a (1,1)-form respectively, by using the method of moving frames.

\vspace{0.2cm}


\begin{center}
\item\section{Introduction}
\end{center}
The twistor construction is an important technique in differential geometry and mathematical physics. In recent years, J. Fine and D. Panov \cite{FP} have introduced the concept of definite connection on $SO(3)$-bundle over an oriented 4-manifold. They showed many non-K\"{a}hler symplectic Calabi-Yau 6-manifolds from the twistor spaces of Riemannian 4-manifolds \cite{FP,FP2,F3}. These works open up a new direction for twistor theory.

The twistor approach was first proposed by R. Penrose in 1960s. In 1978, the Riemannian version of R. Penrose's twistor programme was presented by M. F. Atiyah, N. J. Hitchin and I. M. Singer \cite{AHS}. As showed in \cite{AHS}, basically, to each oriented Riemannian 4-manifold $M$, one can associate canonically a 6-manifold $Z$, the twistor space of $M$, together with an almost complex structure $\mathbb{J}_+$. M. F. Atiyah, N. J. Hitchin and I. M. Singer proved that $\mathbb{J}_+$ is integrable if and only if $M$ is anti-self-dual. Therefore, they established an elegant correspondence between  Yang-Mills fields on 4-manifolds and holomorphic vector bundles on complex 3-manifolds. On the other hand, in 1985, J. Eells and S. Salamon \cite{ES} introduced another almost complex structure $\mathbb{J}_-$ which, by contrast with $\mathbb{J}_+$, is never integrable. However, $\mathbb{J}_-$ plays an important role in the theory of harmonic map.

There is a  natural 1-parameter family of Riemannian metrics $g_t$ on the twistor space $Z$. In fact, $\mathbb{J}_+$ and
$\mathbb{J}_-$ are orthogonal almost complex structures with respect to the metrics $g_t$. Thus, it is natural to study the relations between the almost Hermitian geometry of $Z$ and the Riemannian geometry of 4-manifold $M$ \cite{BN,JR}. For example, N. J. Hitchin \cite{Hi} showed that if the twistor space of  a compact anti-self-dual 4-manifold $M$ admits a K\"{a}hler metric, then $M$ is the 4-sphere $S^4$ or the complex projective plane $\overline{\mathbb{C}P}^2$ (with the non-standard orientation). G. R. Jensen and M. Rigoli \cite{JR} proved that $(Z, g_t, \mathbb{J}_+)$ is an (1,2)-symplectic manifold if and only if $(M, g)$  is  an anti-self-dual, Einstein manifold with  positive scalar curvature $s$ and $t^2=\frac{12}{s}$; while $(Z, g_t, \mathbb{J}_-)$ is an (1,2)-symplectic manifold if and only if $(M, g)$  is  an anti-self-dual, Einstein manifold (for any value of $t>0$).

In this paper, we continue to consider the almost complex structures $\mathbb{J}_+$ and $\mathbb{J}_-$. We use the method of moving frames, which was introduced by G. R. Jensen and M. Rigoli \cite{JR}, to study the almost Hermitian twistor geometry. Some special metric conditions (including Balanced metric condition, first Gauduchon metric condition) on $(Z,g_t,\mathbb{J}_{\pm})$ are studied, which are used to characterize  Riemannian 4-manifold. In particular, we prove that for the Hermitian twistor space $(Z, g_t, \mathbb{J}_+)$ of an anti-self-dual Riemannian 4-manifold $(M, g)$ with constant scalar curvature, $(Z, g_t, \mathbb{J}_+)$ satisfies the first Gauduchon metric condition if and only if $(Z, g_t, \mathbb{J}_+)$ is a K\"{a}hler manifold; and for almost Hermitian twistor space $(Z, g_t, \mathbb{J}_-)$ of an anti-self-dual Riemannian 4-manifold $(M, g)$ with constant scalar curvature, $(Z, g_t, \mathbb{J}_-)$ satisfies the first Gauduchon metric condition if and only if $(M, g)$ is an Einstein manifold.  To some extent, the global version of the first part of this result generalizes a theorem of  M. Verbitsky \cite{Ve}. In fact, our global generalized result is essentially included in \cite{FU,IP}.
Moreover, we also give a natural unitary connection on complex tangent bundle $(TZ,\mathbb{J}_{+})$ with metric $g_t$, and we study the properties of the first Chern forms of the natural complex vector bundles $\mathcal{H}$ and $\mathcal{V}$ over the twistor space $Z$ with the induced unitary connections. some well-known results concerning the first Chern form are recovered by using the method of moving frames.

The paper is organized as follows. In Section 2, we follow the method in \cite{JR} to describe the essential facts of the geometry of Riemannian 4-manifold $M$ in terms of its $SO(4)$-principal bundle of oriented orthonormal frames. In Section 3, we study the almost Hermitian geometry of
the twistor spaces $(Z,g_t,\mathbb{J}_{\pm})$ by using the method of moving frames. In Section 4, unitary connections are constructed on the natural complex vector bundles $\mathcal{H}$ and $\mathcal{V}$ over the twistor space $Z$. We also recover some interesting results with respect to the first Chern form in \cite{FP} and \cite{Gau}.

\begin{center}
\item\section{Preliminaries and notations}
\end{center}
Let $(M,g)$ be an oriented Riemannian 4-manifold. The Hodge star operator gives a map
$\star:\wedge^2\rightarrow \wedge^2$ with $\star^2=1$. Accordingly, its eigenvalues are
$\pm1$ and the bundle of two-forms splits $\wedge^2=\wedge^+\oplus \wedge^-$ into
eigenspaces. $\wedge^+$ (resp. $\wedge^-$) is called the bundle of \emph{self-dual} (resp. \emph{anti-self-dual}) 2-forms.

The Riemannian curvature tensor $R$ of $M$ can be considered as a self-adjoint operator
$\hat{R}:\wedge^2\rightarrow \wedge^2$ and so, with respect to the decomposition $\wedge^2=\wedge^+\oplus \wedge^-$, it decomposes into parts
\begin{equation}
\hat{R}=\left(\begin{array}{cc}A&B^T\\B&C
\end{array}\right),
\end{equation}
where $A=W^{+}+\frac{s}{12}\text{Id}$, $C=W^{-}+\frac{s}{12}\text{Id}$, $W^+$ (resp. $W^-$) is the
self-dual (resp. anti-self-dual) Weyl curvature operator, $s$ is the scalar curvature, and $B=\text{Ric}_0$ is the trace-free Ricci curvature operator \cite{Be}.

Now, we give some fundamental facts of the geometry of Riemannian 4-manifold \cite{JR}, which are useful in the following sections.

The standard action of $SO(4)$ on $\mathbb{R}^4$ (as column vectors) induces a representation
of $SO(4)$ on $\wedge^2\mathbb{R}^4$, which is reducible into irreducible factors $\wedge^2\mathbb{R}^4=\wedge^+\oplus\wedge^-$. For the standard basis $\{\epsilon_1,...,\epsilon_4\}$ of $\mathbb{R}^4$, the standard bases of $\wedge^{\pm}$ are
given by
\begin{equation}
E^{\pm}=(E_1^{\pm},E_2^{\pm},E_3^{\pm}),
\end{equation}
where
$$E_1^{\pm}=\frac{1}{\sqrt{2}}(\epsilon_1\wedge\epsilon_2\pm\epsilon_3\wedge\epsilon_4),
E_2^{\pm}=\frac{1}{\sqrt{2}}(\epsilon_1\wedge\epsilon_3\pm\epsilon_4\wedge\epsilon_2),
E_3^{\pm}=\frac{1}{\sqrt{2}}(\epsilon_1\wedge\epsilon_4\pm\epsilon_2\wedge\epsilon_3).$$

The standard metric on $\mathbb{R}^4$ induces an $SO(4)$-invariant inner product on $\wedge^2\mathbb{R}^4$, and the restriction of the $SO(4)$-action  to $\wedge^{\pm}$ gives a
2:1 surjective homomorphism
\begin{equation}
\iota:SO(4)\rightarrow SO(3)\times SO(3),~~\iota(a)=(a_+,a_-).
\end{equation}

Let $M$ be a connected oriented Riemannian 4-manifold, and  let $P$ denote the $SO(4)$-principal
bundle of oriented orthonormal frames over $M$. We use the index convention $1\leq\alpha,\beta,\gamma,\delta\leq4.$

The $\mathbb{R}^4$-valued canonical form on $P$, denoted by $\theta=(\theta^{\alpha})$, is
given by
$$\theta(X)=e^{-1}(\pi_{\star}(X)),~~~~X\in T_{(x,e)}P,$$
where $e$ is identified with a linear map $e:\mathbb{R}^4\rightarrow T_{x}M$, $\pi:P\rightarrow M$ is the projection.

The $\mathfrak{so}(4)$-valued Levi-Civita connection forms and curvature forms are denoted by
$\omega=(\omega_{\beta}^{\alpha})$ and $\Omega=(\Omega_{\beta}^{\alpha})$, respectively. Then the structure equations of $M$ are
\begin{eqnarray}
&&d\theta^\alpha=-\omega_{\beta}^{\alpha}\wedge\theta^\beta,\\
&&d\omega_{\beta}^{\alpha}=-\omega_{\gamma}^{\alpha}\wedge\omega_{\beta}^{\gamma}+\Omega_{\beta}^{\alpha},
\end{eqnarray}
where $\Omega_{\beta}^{\alpha}=\frac{1}{2}R_{\alpha\beta\gamma\delta}\theta^\gamma\wedge\theta^\delta$,
$R_{\alpha\beta\gamma\delta}$ are functions on $P$ defining the Riemannian curvature tensor of $M$.

For any $a\in SO(4)$, we have \cite{KN}
\begin{eqnarray}
&&R_a^{\star}\theta=a^{-1}\theta,\\
&&R_a^{\star}\omega=Ad(a^{-1})\omega=a^{-1}\omega a,\\
&&R_a^{\star}\Omega=Ad(a^{-1})\Omega=a^{-1}\Omega a,
\end{eqnarray}
where $R_a$ denotes the right multiplication on $P$.

We define $\mathbb{R}^3$-valued 2-forms on $P$ by
\begin{equation}
\alpha_{\pm}=(\alpha_{\pm}^1,\alpha_{\pm}^2,\alpha_{\pm}^3)^T,
\end{equation}
where
$$\alpha_{\pm}^1=\frac{1}{\sqrt{2}}(\theta^1\wedge\theta^2\pm\theta^3\wedge\theta^4),
\alpha_{\pm}^2=\frac{1}{\sqrt{2}}(\theta^1\wedge\theta^3\pm\theta^4\wedge\theta^2),
\alpha_{\pm}^3=\frac{1}{\sqrt{2}}(\theta^1\wedge\theta^4\pm\theta^2\wedge\theta^3).$$

For the curvature forms matrix $\Omega$, in terms of the basis (2.2) and (2.9),
$$\Omega=E^+\otimes A \alpha_++E^-\otimes B \alpha_+
+E^+\otimes B^T \alpha_-+E^-\otimes C \alpha_-,$$
where $A=(A_{ij}),~B=(B_{ij})$ and $C=(C_{ij})$ are $3\times 3$-matrix-valued functions on $P$, and $A^T=A$, $C^T=C$. For the explicit formulas relating $R_{\alpha\beta\gamma\delta}$ to $A,~B,~C$, refer to (3.11)-(3.13) in \cite{JR}. For any $a\in SO(4)$, it follows from (2.8) that \cite{JR}
\begin{equation}
R_a^{\star} A=a_+^{-1}A a_+,~R_a^{\star} B=a_-^{-1}B a_+,~R_a^{\star} C=a_-^{-1}C a_-.
\end{equation}

It is well-known that $(M,g)$ is an Einstein manifold if and only if $B=0$ on $P$. An oriented Riemannian 4-manifold $(M,g)$ is called \emph{anti-self-dual} if the self-dual Weyl curvature operator $W^+=0$ on $M$. In fact, $(M,g)$ is an anti-self-dual manifold is equivalent to
$A-\frac{s}{12}\text{Id}=0$ on $P$.

Set $Z=\{(x,J_{x})|x\in M, J_{x}$ is an
orientation preserving orthogonal complex structure of the vector space $T_{x}M\}$.
$Z$ is called the \emph{twistor space} \cite{AHS} of $(M,g)$. $Z$ has the following representations:
$$Z=P\times_{SO(4)}SO(4)/U(2)=S(\wedge^{+}),$$
where $S(\wedge^{+})$ is the sphere bundle associated to the bundle of self-dual 2-forms $\wedge^{+}$.

From the definition of twistor space, it follows that $Z$ depends only on the conformal
structure of $M$. Notations: $\pi:P\rightarrow M$, $\pi_1:P\rightarrow Z$ and $\pi_2:Z\rightarrow M$ are the
projections.

We first consider the representation of $U(2)$ in $SO(4)$. Let $E_{12}$ be an $4\times4$-matrix
with (1,2)-component -1, (2,1)-component 1, other components are
zeros,  and let $E_{34}$ be an $4\times4$-matrix
with (3,4)-component -1, (4,3)-component 1, other components are
zeros. Set $J_2=E_{12}+E_{34}$, then
\begin{equation}
\varrho:U(2)\cong\{a\in SO(4)|a^TJ_2a=J_2\}.
\end{equation}
The corresponding Lie algebra representation is
$$\mathfrak{u}(2)\cong\{X\in \mathfrak{so}(4)|X^TJ_2+J_2X=0\}.$$

Up to a positive factor, the unique $Ad(SO(4))$-invariant inner product on $\mathfrak{so}(4)$ is
$$<X,Y>=\text{trace}~(X^T Y),$$
where $X,Y\in\mathfrak{so}(4)$.

Let $\mathfrak{m}$ be the orthogonal complement of
$\mathfrak{u}(2)$ in $\mathfrak{so}(4)$. Then $\mathfrak{so}(4)=\mathfrak{u}(2)\oplus\mathfrak{m}$, and the connection forms matrix
$\omega=\mu+\nu$, where
$$\mu=\frac{1}{2}(\omega-J_2\omega J_2),~~\nu=\frac{1}{2}(\omega+J_2\omega J_2).$$

We define a symmetric bilinear form $P_t$ on $P$ by
\begin{equation}
P_t=\theta^T\cdot\theta+4t^2((\theta^5)^2+(\theta^6)^2),
\end{equation}
where $t>0$, $\theta^5=\frac{1}{2}(\omega_3^1-\omega_4^2),\theta^6=\frac{1}{2}(\omega_4^1+\omega_3^2).$

It is easily checked that $R_{\varrho(a)}^{\star}P_t=P_t$ for any $a\in U(2)$, and that $P_t$ is horizontal. Thus there exists a unique Riemannian metric $g_t$ on $Z$ such that $\pi_1^\star g_t=P_t$ \cite{JR}.

Let $U\subset Z$ be an open subset on which there is a local section $u:U\rightarrow P$. An
orthonormal co-frame for $g_t$ on $U$ is given  by
\begin{equation}
u^{\star}(\theta^\alpha),~2tu^{\star}(\theta^5),~2tu^{\star}(\theta^6).
\end{equation}
With respect to this orthonormal co-frame, we can calculate the components of the Riemannian curvature tensor and study the geometry of the twistor space $Z$. In the following section, we consider some special metric conditions on the almost Hermitian twistor space by using the method of moving frames introduced in \cite{JR}.

\begin{center}
\item\section{Almost Hermitian geometry of the twistor space}
\end{center}
The basic fact of twistor theory is that the twistor space $Z$ has two
distinguished almost complex structures defined by the Levi-Civita connection of
Riemannian 4-manifold $M$, denoted by $\mathbb{J}_{\pm}$, which are introduced by
Atiyah-Hitchin-Singer and Eells-Salamon, respectively.

G. R. Jensen and M. Rigoli \cite{JR}  also described  almost complex structures $\mathbb{J}_{\pm}$ on $Z$ and studied their almost Hermitian twistor geometry by the language of principal bundle.

As in \cite{JR}, to describe $\mathbb{J}_{\pm}$, we do this by defining locally (1,0)-forms on $Z$.
On the $SO(4)$-principal bundle $P$, we have canonical form $\theta$ and the Levi-Civita connection form $\omega$. We define complex-valued 1-forms on $P$ by
$$\varphi^1=\theta^1+\sqrt{-1}\theta^2,~\varphi^2=\theta^3+\sqrt{-1}\theta^4,~\varphi^3
=\theta^5+\sqrt{-1}\theta^6.$$

Let $U\subset Z$ be an open subset on which there is a local section $u:U\rightarrow P$. It follows  that the local complex-valued 1-forms
$\{u^{\star}\varphi^1$,  $u^{\star}\varphi^2$,  $u^{\star}\varphi^3\}$ are the basis of (1,0)-forms of the almost complex structure $\mathbb{J}_+$, and the local complex-valued 1-forms
$\{u^{\star}\varphi^1$,  $u^{\star}\varphi^2$,  $u^{\star}\bar{\varphi}^3\}$ are the basis of (1,0)-forms of the almost complex structure $\mathbb{J}_-$. In the following, $(Z,g_t,\mathbb{J}_{+})$ and $(Z,g_t,\mathbb{J}_{-})$ (or abbr. $(Z,g_t,\mathbb{J}_{\pm})$) are called the \emph{almost Hermitian twistor spaces}.

By the structure equations (2.4) and (2.5), on the principal bundle $P$, we have
\begin{eqnarray}
&&d\varphi^3=\sqrt{-1}(\omega_2^1+\omega_4^3)\wedge\varphi^3+\rho(\Omega),\\
&&d\bar{\varphi}^3=-\sqrt{-1}(\omega_2^1+\omega_4^3)\wedge\bar{\varphi}^3
+\bar{\rho}(\Omega),\\
&&d\varphi^1=\sqrt{-1}\omega_2^1\wedge\varphi^1-\psi^3\wedge\varphi^2
+\bar{\varphi}^2\wedge\varphi^3,\\
&&d\varphi^2=\bar{\psi}^3\wedge\varphi^1+\sqrt{-1}\omega_4^3\wedge\varphi^2
+\varphi^3\wedge\bar{\varphi}^1,
\end{eqnarray}
where
$$\psi^3=\frac{1}{2}[\omega_3^1+\omega_4^2+\sqrt{-1}(\omega_3^2-\omega_4^1)],$$
$$\rho(\Omega)=\frac{1}{2}[\Omega_3^1-\Omega_4^2+\sqrt{-1}(\Omega_3^2+\Omega_4^1)].$$

For the almost Hermitian twistor spaces $(Z,g_t,\mathbb{J}_{\pm})$, the K\"{a}hler forms, denoted
by $\mathbb{K}_{\pm}(t)$, are
\begin{equation}
\mathbb{K}_{\pm}(t)=\frac{\sqrt{-1}}{2}(u^{\star}\varphi^1\wedge u^{\star}\bar{\varphi}^1+
u^{\star}\varphi^2\wedge u^{\star}\bar{\varphi}^2\pm4t^2u^{\star}\varphi^3\wedge u^{\star}\bar{\varphi}^3),
\end{equation}
where $u$ is any local section of $\pi_1:P\rightarrow Z$.

Hereafter, for convenience, we always omit the pullback mapping $u^{\star}$ and leave it to the context to show the manifold on which the forms are considered.

Set
\begin{equation}
\rho(\Omega)=\lambda_{12}\varphi^1\wedge\varphi^2+
\lambda_{\bar{1}\bar{2}}\bar{\varphi}^1\wedge\bar{\varphi}^2+\lambda_{1\bar{1}}\varphi^1\wedge\bar{\varphi}^1
\end{equation}
~~~~~~~~~~~~~~~~~~~~~~~~~~~~~$+\lambda_{2\bar{2}}\varphi^2\wedge\bar{\varphi}^2
+\lambda_{\bar{1}2}\bar{\varphi}^1\wedge\varphi^2+
\lambda_{1\bar{2}}\varphi^1\wedge\bar{\varphi}^2.
$

\noindent Since $\Omega_{\beta}^{\alpha}=\frac{1}{2}R_{\alpha\beta\gamma\delta}\theta^\gamma\wedge\theta^\delta$, it follows
$$\lambda_{12}=\frac{1}{4}(A_{22}+A_{33}),~\lambda_{\bar{1}\bar{2}}=\frac{1}{4}(A_{22}-A_{33}
+2\sqrt{-1}A_{23}),$$
\begin{equation}
\lambda_{1\bar{1}}=\frac{1}{4}(\sqrt{-1}(A_{12}+B_{12})-A_{13}-B_{13}),~
\lambda_{2\bar{2}}=\frac{1}{4}(\sqrt{-1}(A_{12}-B_{12})-A_{13}+B_{13}),
\end{equation}
$$
\lambda_{\bar{1}2}=\frac{1}{4}(\sqrt{-1}(B_{23}-B_{32})+B_{22}+B_{33}),~
\lambda_{1\bar{2}}=\frac{1}{4}(\sqrt{-1}(B_{23}+B_{32})+B_{22}-B_{33}).
$$
\noindent Then we have
\begin{equation}
d\mathbb{K}_{\pm}(t)=-\sqrt{-1}\varphi^3\wedge[(-1\pm2t^2\lambda_{12})\bar{\varphi}^1\wedge\bar{\varphi}^2
\pm2t^2\bar{\Lambda}]
\end{equation}
~~~~~~~~~~~~~~~~~~~~~~~~~~~~~$+\sqrt{-1}\bar{\varphi}^3\wedge[(-1\pm2t^2\lambda_{12})\varphi^1\wedge\varphi^2
\pm2t^2\Lambda],$

\noindent where $\Lambda=\rho(\Omega)-\lambda_{12}\varphi^1\wedge\varphi^2.$

The (1,2)-symplectic condition on $(Z,g_t,\mathbb{J}_{\pm})$ is studied by G. R. Jensen and M. Rigoli \cite{JR}. In fact, if $(Z,g_t,\mathbb{J}_+)$ (resp. $(Z,g_t,\mathbb{J}_-)$) is an (1,2)-symplectic manifold, then $\mathbb{K}_{+}(t)\wedge d\mathbb{K}_{+}(t)=0$ (resp. $\mathbb{K}_{-}(t)\wedge d\mathbb{K}_{-}(t)=0$). Conversely, we have the following theorem:\vspace{0.2cm}

\noindent\textbf{Theorem 3.1} ~~\emph{Let $(M, g)$  be an oriented Riemannian 4-manifold, its almost Hermitian twistor spaces are $(Z,g_t,\mathbb{J}_{\pm})$, the K\"{a}hler forms are $\mathbb{K}_{\pm}(t)$.  The following conditions are equivalent:}

\emph{(I)~~~$(Z, g_t, \mathbb{J}_{+})$ satisfies $\mathbb{K}_{+}(t)\wedge d\mathbb{K}_{+}(t)=0$;}

\emph{(II)~~$(Z, g_t, \mathbb{J}_{-})$ satisfies $\mathbb{K}_{-}(t)\wedge d\mathbb{K}_{-}(t)=0$;}

\emph{(III)~$(M, g)$  is an anti-self-dual manifold.} \vspace{0.2cm}

\emph{Proof.} From (3.5) and (3.8), we have
$$\mathbb{K}_{+}(t)\wedge d\mathbb{K}_{+}(t)=t^2\{(\bar{\lambda}_{1\bar{1}}+\bar{\lambda}_{2\bar{2}})\varphi^1\wedge\bar{\varphi}^1
\wedge\bar{\varphi}^2\wedge\varphi^2\wedge\varphi^3$$
~~~~~~~~~~~~~~~~~~~~~~~~~~~~~~~~~~~~~$+(\lambda_{1\bar{1}}+\lambda_{2\bar{2}})\varphi^1\wedge\bar{\varphi}^1
\wedge\bar{\varphi}^2\wedge\varphi^2\wedge\bar{\varphi}^3\},$

$$\mathbb{K}_{-}(t)\wedge d\mathbb{K}_{-}(t)=-t^2\{(\bar{\lambda}_{1\bar{1}}+\bar{\lambda}_{2\bar{2}})\varphi^1\wedge\bar{\varphi}^1
\wedge\bar{\varphi}^2\wedge\varphi^2\wedge\varphi^3$$
~~~~~~~~~~~~~~~~~~~~~~~~~~~~~~~~~~~~~$+(\lambda_{1\bar{1}}+\lambda_{2\bar{2}})\varphi^1\wedge\bar{\varphi}^1
\wedge\bar{\varphi}^2\wedge\varphi^2\wedge\bar{\varphi}^3\}.$\vspace{0.2cm}

From the above two equalities, it follows that $\mathbb{K}_{+}(t)\wedge d\mathbb{K}_{+}(t)=-\mathbb{K}_{-}(t)\wedge d\mathbb{K}_{-}(t)$, i.e., (I) is
equivalent to (II).

$\mathbb{K}_{+}(t)\wedge d\mathbb{K}_{+}(t)=0$ if and only if $\lambda_{1\bar{1}}+\lambda_{2\bar{2}}=0$ on $P$. Using (3.7), we obtain $\mathbb{K}_{+}(t)\wedge d\mathbb{K}_{+}(t)=0$ if and only if $A_{12}=A_{13}=0$ on $P$.

The method is the same as in the proof of Theorem 6.1 in \cite{JR}. By the transformation equalities (2.10), and the fact that that the homomorphism  $\iota:SO(4)\rightarrow SO(3)\times SO(3),~\iota(a)=(a_+,a_-)$ is surjective, it follows that
$A_{12}=A_{13}=0$ on $P$ if and only if $A$ is a scalar matrix.

As $\text{trace}~A=\frac{s}{4}$, it follows that $A$ is a scalar matrix if and only if $A=\frac{s}{12}\text{Id}$ on $P$, i.e., $(M, g)$  is an anti-self-dual manifold.
\hfill$\Box$
\vspace{0.2cm}

\noindent\textbf{Remark 3.2}~~The metric condition on $(Z,g_t,\mathbb{J}_{\pm})$ in Theorem 3.1 is called the Balanced metric condition. M. L. Michelsohn \cite{Mi} introduced the Balanced metric on complex manifold, and also claimed that the natural metric on the twistor space of an anti-self-dual 4-manifold is a Balanced metric.\vspace{0.2cm}

$(Z,g_t,\mathbb{J}_+)$ is a Hermitian manifold only when $(M, g)$  is an anti-self-dual 4-manifold, while $(Z,g_t,\mathbb{J}_-)$ is not a  Hermitian manifold in any case. We want to
consider more special metric conditions on the twistor space, which enable us to obtain some metric properties on $(M, g)$.

We begin with simple remarks concerning the decomposition of the exterior differentiation $d$ in the almost complex case \cite{Br}. In general, if $\phi$ is a $(p,q)$-form, then $d\phi$ decomposes into
a sum of forms of degree $(p+2,q-1)$, $(p+1,q)$, $(p,q+1)$, $(p-1,q+2)$. Thus $d\phi=d^{2,-1}\phi+\partial\phi+\bar{\partial}\phi+d^{-1,2}\phi$.

For almost Hermitian twistor space $(Z,g_t,\mathbb{J}_+)$, form (3.8), we have
\begin{eqnarray*}
&&d^{2,-1}\mathbb{K}_{+}(t)=-2\sqrt{-1}t^2\bar{\lambda}_{\bar{1}\bar{2}}\varphi^1\wedge\varphi^2
\wedge\varphi^3,\\
&&d^{-1,2}\mathbb{K}_{+}(t)=2\sqrt{-1}t^2\lambda_{\bar{1}\bar{2}}\bar{\varphi}^1\wedge\bar{\varphi}^2
\wedge\bar{\varphi}^3,\\
&&\partial\mathbb{K}_{+}(t)=\sqrt{-1}(2t^2\lambda_{12}-1)\varphi^1\wedge\varphi^2
\wedge\bar{\varphi}^3-2\sqrt{-1}t^2\varphi^3\wedge[\bar{\Lambda}
-\bar{\lambda}_{\bar{1}\bar{2}}\varphi^1\wedge\varphi^2],\\
&&\bar{\partial}\mathbb{K}_{+}(t)=\sqrt{-1}(1-2t^2\lambda_{12})\bar{\varphi}^1\wedge\bar{\varphi}^2
\wedge\varphi^3+2\sqrt{-1}t^2\bar{\varphi}^3\wedge[\Lambda
-\lambda_{\bar{1}\bar{2}}\bar{\varphi}^1\wedge\bar{\varphi}^2].
\end{eqnarray*}

For almost Hermitian twistor space $(Z,g_t,\mathbb{J}_-)$, from (3.8), we have
\begin{eqnarray*}
&&d^{2,-1}\mathbb{K}_{-}(t)=-\sqrt{-1}(1+2t^2\lambda_{12})\varphi^1\wedge\varphi^2
\wedge\bar{\varphi}^3,\\
&&d^{-1,2}\mathbb{K}_{-}(t)=\sqrt{-1}(1+2t^2\lambda_{12})\bar{\varphi}^1\wedge\bar{\varphi}^2
\wedge\varphi^3,\\
&&\partial\mathbb{K}_{-}(t)=2\sqrt{-1}t^2\bar{\lambda}_{\bar{1}\bar{2}}\varphi^1\wedge\varphi^2
\wedge\varphi^3-2\sqrt{-1}t^2\bar{\varphi}^3\wedge[\Lambda
-\lambda_{\bar{1}\bar{2}}\bar{\varphi}^1\wedge\bar{\varphi}^2],\\
&&\bar{\partial}\mathbb{K}_{-}(t)=2\sqrt{-1}t^2\varphi^3\wedge[\bar{\Lambda}
-\bar{\lambda}_{\bar{1}\bar{2}}\varphi^1\wedge\varphi^2]-
2\sqrt{-1}t^2\lambda_{\bar{1}\bar{2}}\bar{\varphi}^1\wedge\bar{\varphi}^2
\wedge\bar{\varphi}^3.
\end{eqnarray*}

Since $d^{2,-1}$ and $d^{-1,2}$ are linear over the smooth functions, from the structure equations (3.1)-(3.4), and by direct calculations, we get the following proposition:\vspace{0.2cm}

\noindent\textbf{Proposition 3.3}~~\emph{For almost Hermitian twistor spaces $(Z,g_t,\mathbb{J}_{\pm})$, the K\"{a}hler forms are $\mathbb{K}_{\pm}(t)$, then}
\begin{eqnarray*}
&&d^{2,-1}\bar{\partial}\mathbb{K}_{+}(t)=0,\\
&&d^{2,-1}d^{-1,2}\mathbb{K}_{+}(t)=2\sqrt{-1}t^2|\lambda_{\bar{1}\bar{2}}|^2
\varphi^1\wedge\varphi^2\wedge\bar{\varphi}^1\wedge\bar{\varphi}^2,\\
&&d^{2,-1}\bar{\partial}\mathbb{K}_{-}(t)=2\sqrt{-1}t^2(\bar{\lambda}_{1\bar{1}}
+\bar{\lambda}_{2\bar{2}})\varphi^1\wedge\varphi^2\wedge\bar{\varphi}^3\wedge\varphi^3,\\
&&d^{2,-1}d^{-1,2}\mathbb{K}_{-}(t)=\sqrt{-1}(1+2t^2\lambda_{12})(\lambda_{12}
\bar{\varphi}^1\wedge\varphi^1\wedge\varphi^2\wedge\bar{\varphi}^2
\end{eqnarray*}
~~~~~~~~~~~~~~~~~~~~~~~~~~~~~~~~~~~~~~~~$+\bar{\varphi}^2\wedge\varphi^2\wedge\bar{\varphi}^3
\wedge\varphi^3
+\varphi^3\wedge\bar{\varphi}^3\wedge\varphi^1\wedge\bar{\varphi}^1).$\vspace{0.2cm}

\noindent\textbf{Remark 3.4}~~For the general case, the formulae of $\partial\bar{\partial}\mathbb{K}_{\pm}(t)$ are complicated.
For anti-self-dual 4-manifold, A. Nannicini \cite{Na} showed the formula of  $\partial\bar{\partial}\mathbb{K}_{+}(t)$ on local  horizontal frame and
local vertical frame. For anti-self-dual 4-manifold with constant scalar curvature,  G. Deschamps, N. Le Du and C. Mourougane \cite{DLM} also compute  $\partial\bar{\partial}\mathbb{K}_{+}(t)$. Meanwhile, they  analyse the condition  $\partial\bar{\partial}\mathbb{K}_{+}(t)=0$, and recover a local version of the result of M. Verbitsky \cite{Ve}.
\vspace{0.2cm}

Given a complex $n$-dimensional Hermitian manifold with metric $ds^2$, its K\"{a}hler form is denoted by $\mathbb{K}$. The  metric $ds^2$ is called a \emph{pluriclosed metric} or \emph{strong K\"{a}hler with torsion metric} if $\partial\bar{\partial}\mathbb{K}=0$. The  metric $ds^2$ is called a \emph{Gauduchon  metric}
if $\partial\bar{\partial}\mathbb{K}^{n-1}=0$. Z. Z. Wang, D. M. Wu and the first named author \cite{FWW} introduced the generalized Gauduchon metric, that is, $ds^2$ is called a \emph{$l$-th Gauduchon metric} if
$\partial\bar{\partial}\mathbb{K}^l\wedge \mathbb{K}^{n-l-1}=0$, $1\leq l\leq n-1$.

For an anti-self-dual 4-manifold $M$, on the twistor space $Z$, from Theorem 3.1, we have
$\mathbb{K}_{+}(t)\wedge\bar{\partial}\mathbb{K}_{+}(t)=0$, $\mathbb{K}_{-}(t)\wedge\bar{\partial}\mathbb{K}_{-}(t)=0$. Thus, it follows that
$\partial\bar{\partial}\mathbb{K}_{\pm}^{2}=0$, i.e., the natural metrics $g_t$ on $(Z,\mathbb{J}_{\pm})$ are Gauduchon metrics.

In the following, we study the first Gauduchon metrics (i.e., the condition
$\partial\bar{\partial}\mathbb{K}_{\pm}(t)\wedge \mathbb{K}_{\pm}(t)=0$) on the twistor spaces $(Z,g_t,\mathbb{J}_{\pm})$. In particular, we generalize a result in \cite{Ve} and \cite{DLM} by using a simple method.

Here, we abuse the definition of $l$-th Gauduchon metric on the almost Hermitian manifold as on the Hermitian case.\vspace{0.2cm}

\noindent\textbf{Theorem 3.5}~~\emph{Let $(M, g)$  be an anti-self-dual Riemannian 4-manifold with
constant scalar curvature $s$,  the corresponding Hermitian twistor space is $(Z,g_t,\mathbb{J}_{+})$, the K\"{a}hler form is denoted by $\mathbb{K}_{+}(t)$. Then
$\partial\bar{\partial}\mathbb{K}_{+}(t)\wedge \mathbb{K}_{+}(t)=0$ if and only if
$d\mathbb{K}_{+}(t)=0$, i.e., $(Z,g_t,\mathbb{J}_{+})$ is a K\"{a}hler manifold.}
\vspace{0.2cm}

\emph{Proof.} From the proof of Theorem 5.23 in \cite{JR}, $(M, g)$ is an anti-self-dual
 manifold if and only if
$\lambda_{\bar{1}\bar{2}}=0$ on principal bundle $P$. Under the hypothesis condition, it follows
$$\bar{\partial}\mathbb{K}_{+}(t)=\sqrt{-1}(1-2t^2\lambda_{12})\bar{\varphi}^1\wedge\bar{\varphi}^2
\wedge\varphi^3+2\sqrt{-1}t^2\bar{\varphi}^3\wedge(\rho(\Omega)
-\lambda_{12}\varphi^1\wedge\varphi^2),$$
where $\lambda_{12}=\frac{s}{24}$ is a constant.

By the structure equations  (3.1)-(3.4), we have
\begin{equation}
d(\bar{\varphi}^1\wedge\bar{\varphi}^2
\wedge\varphi^3)=\varphi^2\wedge\bar{\varphi}^2\wedge\varphi^3\wedge\bar{\varphi}^3
+\varphi^1\wedge\bar{\varphi}^1\wedge\varphi^3\wedge\bar{\varphi}^3
+\bar{\varphi}^1\wedge\bar{\varphi}^2\wedge \rho(\Omega).
\end{equation}
From (2.5) and its exterior differentiation,  we obtain
\begin{equation}
d\rho(\Omega)=\sqrt{-1}(\omega_2^1+\omega_4^3)\wedge \rho(\Omega)-\sqrt{-1}(\Omega_2^1+\Omega_4^3)\wedge \varphi^3.
\end{equation}

Thus, together with (3.2), (3.9) and (3.10), we have
\begin{equation}
\partial\bar{\partial}\mathbb{K}_{+}(t)=
d\bar{\partial}\mathbb{K}_{+}(t)=\sqrt{-1}(1-4t^2\lambda_{12})\varphi^2\wedge\bar{\varphi}^2\wedge\varphi^3\wedge \bar{\varphi}^3
\end{equation}

~~~~~~~~~~~~~~~~~~~~~~~~~~~~~~~~~$+\sqrt{-1}(1-4t^2\lambda_{12})
\varphi^1\wedge\bar{\varphi}^1\wedge\varphi^3\wedge\bar{\varphi}^3$

~~~~~~~~~~~~~~~~~~~~~~~~~~~~~~~~~$+\sqrt{-1}\lambda_{12}(1-4t^2\lambda_{12})
\varphi^1\wedge\varphi^2\wedge\bar{\varphi}^1\wedge\bar{\varphi}^2$

~~~~~~~~~~~~~~~~~~~~~~~~~~~~~~~~~$+2t^2(\Omega_2^1+\Omega_4^3)\wedge \varphi^3\wedge \bar{\varphi}^3+2\sqrt{-1}t^2 \rho(\Omega)\wedge \bar{\rho}(\Omega).$

\noindent From (3.6), it follows
\begin{equation}
\rho(\Omega)\wedge \bar{\rho}(\Omega)=(\lambda_{12}^2-|\lambda_{1\bar{1}}|^2
-|\lambda_{2\bar{2}}|^2-|\lambda_{\bar{1}2}|^2-|\lambda_{1\bar{2}}|^2
)\varphi^1\wedge\varphi^2\wedge\bar{\varphi}^1\wedge\bar{\varphi}^2,
\end{equation}
then, together with (3.11), we obtain
\begin{equation}
\partial\bar{\partial}\mathbb{K}_{+}(t)=\sqrt{-1}[1-4t^2\lambda_{12}+t^2(R_{3434}+R_{1234})]
\varphi^2\wedge\bar{\varphi}^2\wedge\varphi^3\wedge \bar{\varphi}^3
\end{equation}

~~~~~~~~$+\sqrt{-1}[1-4t^2\lambda_{12}+t^2(R_{3412}+R_{1212})]
\varphi^1\wedge\bar{\varphi}^1\wedge\varphi^3\wedge \bar{\varphi}^3$

~~~~~~~~~~$+\sqrt{-1}\Lambda_1\varphi^1\wedge\varphi^2\wedge\bar{\varphi}^1\wedge\bar{\varphi}^2
+\Lambda_2\varphi^1\wedge\bar{\varphi}^2\wedge\varphi^3\wedge \bar{\varphi}^3
+\bar{\Lambda}_2\bar{\varphi}^1\wedge\varphi^2\wedge\varphi^3\wedge \bar{\varphi}^3,$

\noindent where
$$\Lambda_1=\lambda_{12}(1-2t^2\lambda_{12})-2t^2(|\lambda_{1\bar{1}}|^2
+|\lambda_{2\bar{2}}|^2+|\lambda_{\bar{1}2}|^2+|\lambda_{1\bar{2}}|^2),$$
$$\Lambda_2=\frac{t^2}{2}[R_{1213}+R_{3413}+R_{1224}+R_{3424}
+\sqrt{-1}(R_{1214}+R_{3414}-R_{1223}-R_{3423})],$$

It follows that $\partial\bar{\partial}\mathbb{K}_{+}(t)\wedge \mathbb{K}_{+}(t)=0$ if and only if \begin{equation}
2-8t^2\lambda_{12}+t^2(R_{1212}+R_{3434}+2R_{1234})-4t^2\Lambda_1=0
\end{equation}
on principal bundle $P$.

As $M$ is anti-self-dual, then $A=\frac{s}{12}\text{Id}$, $2A_{11}=R_{1212}+R_{3434}+2R_{1234}$ and $2\lambda_{12}=A_{11}$. Through simplification, equality (3.14) transforms to
\begin{equation}
2(1-t^2A_{11})^2+8t^4(|\lambda_{1\bar{1}}|^2
+|\lambda_{2\bar{2}}|^2+|\lambda_{\bar{1}2}|^2+|\lambda_{1\bar{2}}|^2)=0.
\end{equation}
Thus, $\partial\bar{\partial}\mathbb{K}_{+}(t)\wedge \mathbb{K}_{+}(t)=0$ if and only if
\begin{equation}
1-t^2A_{11}=1-2t^2\lambda_{12}=0,~~\lambda_{1\bar{1}}=\lambda_{2\bar{2}}=\lambda_{\bar{1}2}=
\lambda_{1\bar{2}}=0
\end{equation}
on principal bundle $P$.

From (3.8) and (3.16), it follows that
$\partial\bar{\partial}\mathbb{K}_{+}(t)\wedge \mathbb{K}_{+}(t)=0$ if and only if
$d\mathbb{K}_{+}(t)=0$.
\hfill$\Box$
\vspace{0.2cm}

From Theorem 3.5 and  N. J. Hitchin's \cite{Hi} classification theorem of K\"{a}hler twistor spaces, we obtain the following global result which is essentially included in Proposition 2.4 in \cite{FU} or Theorem 1.3 in \cite{IP}.
\vspace{0.2cm}

\noindent\textbf{Corollary 3.6}~~\emph{Let $(M, g)$  be a compact anti-self-dual Riemannian 4-manifold with constant scalar curvature $s$,  the corresponding Hermitian twistor space is $(Z,g_t,\mathbb{J}_{+})$, the K\"{a}hler form is denoted by $\mathbb{K}_{+}(t)$. If
$\partial\bar{\partial}\mathbb{K}_{+}(t)\wedge \mathbb{K}_{+}(t)=0$, then $M$ is isomorphic to
the complex projective plane $\overline{\mathbb{C}P}^2$ (with the
non-standard orientation) or the 4-sphere $S^4$.}
\vspace{0.2cm}

For the  almost Hermitian twistor space  $(Z,g_t,\mathbb{J}_{-})$, we have
\vspace{0.2cm}

\noindent\textbf{Theorem 3.7}~~\emph{Let $(M, g)$  be an anti-self-dual Riemannian 4-manifold with
constant scalar curvature $s$,  the corresponding almost Hermitian twistor space is $(Z,g_t,\mathbb{J}_{-})$, the K\"{a}hler form is denoted by $\mathbb{K}_{-}(t)$. Then
$\partial\bar{\partial}\mathbb{K}_{-}(t)\wedge \mathbb{K}_{-}(t)=0$ if and only if
$(M, g)$ is an Einstein manifold.}
\vspace{0.2cm}

\emph{Proof.} The method is the same as the proof of Theorem 3.5.  For completeness, we give a brief proof.
As $M$ is anti-self-dual, then $\lambda_{\bar{1}\bar{2}}=0$ on principal bundle $P$ and
$$\bar{\partial}\mathbb{K}_{-}(t)=2\sqrt{-1}t^2[\bar{\rho}(\Omega)\wedge \varphi^3
-\lambda_{12}\bar{\varphi}^1\wedge\bar{\varphi}^2
\wedge\varphi^3],$$
where $\lambda_{12}=\frac{s}{24}$ is a constant.

From (3.10) and the structure equations (3.1)-(3.4), we have
\begin{equation}
d\bar{\partial}\mathbb{K}_{-}(t)=2t^2(\Omega_2^1+\Omega_4^3)\wedge\varphi^3\wedge\bar{\varphi}^3
-2\sqrt{-1}t^2\lambda_{12}(\varphi^2\wedge\bar{\varphi}^2\wedge\varphi^3\wedge\bar{\varphi}^3
+\varphi^1\wedge\bar{\varphi}^1\wedge\varphi^3\wedge\bar{\varphi}^3)
\end{equation}

~~~~~~~~~~~$+2\sqrt{-1}t^2(|\lambda_{1\bar{1}}|^2
+|\lambda_{2\bar{2}}|^2+|\lambda_{\bar{1}2}|^2+|\lambda_{1\bar{2}}|^2)
\varphi^1\wedge\bar{\varphi}^1\wedge\varphi^2\wedge\bar{\varphi}^2$.

In fact, as $M$ is anti-self-dual, it follows
\begin{equation}
\Omega_2^1+\Omega_4^3=\frac{\sqrt{-1}}{2}[(R_{1212}+R_{3412})\varphi^1\wedge\bar{\varphi}^1
+(R_{3434}+R_{1234})\varphi^2\wedge\bar{\varphi}^2]
\end{equation}

~~~~~~~~~~~~~~$+\frac{1}{2t^2}(\Lambda_2\varphi^1\wedge\bar{\varphi}^2
+\bar{\Lambda}_2\bar{\varphi}^1\wedge\varphi^2)$.

Then, for almost Hermitian structure $\mathbb{J}_-$, from (3.17) and (3.18), we also have
$$\partial\bar{\partial}\mathbb{K}_{-}(t)=d\bar{\partial}\mathbb{K}_{-}(t).$$
By direct calculations as in the proof of Theorem 3.5, it follows that
$\partial\bar{\partial}\mathbb{K}_{-}(t)\wedge \mathbb{K}_{-}(t)=0$ if and only if
\begin{equation}
-4t^2\lambda_{12}+t^2(R_{1212}+R_{3434}+2R_{1234})-8t^4(|\lambda_{1\bar{1}}|^2
+|\lambda_{2\bar{2}}|^2+|\lambda_{\bar{1}2}|^2+|\lambda_{1\bar{2}}|^2)=0
\end{equation}
on principal bundle $P$.

For anti-self-dual manifold, $4\lambda_{12}=R_{1212}+R_{3434}+2R_{1234}$, thus $\partial\bar{\partial}\mathbb{K}_{-}(t)\wedge \mathbb{K}_{-}(t)=0$ if and only if
\begin{equation}
\lambda_{1\bar{1}}=\lambda_{2\bar{2}}=\lambda_{\bar{1}2}=\lambda_{1\bar{2}}=0
\end{equation}
on principal bundle $P$.

From the proof of Theorem 6.1 in \cite{JR}, for anti-self-dual manifold, $\lambda_{1\bar{1}}=\lambda_{2\bar{2}}=\lambda_{\bar{1}2}=\lambda_{1\bar{2}}=0$
on principal bundle $P$ if and only if $B=0$, i.e., $M$ is an Einstein manifold.
\hfill$\Box$
\vspace{0.2cm}

\noindent\textbf{Remark 3.8}~~For an anti-self-dual Einstein 4-manifold $(M, g)$, the corresponding K\"{a}hler forms
$\mathbb{K}_{\pm}(t)$ of the almost Hermitian twistor spaces $(Z,g_t,\mathbb{J}_{\pm})$ satisfy the following equalities:\vspace{0.2cm}

$\bar{\partial}\mathbb{K}_{-}(t)=\partial\mathbb{K}_{-}(t)=0,$

$\partial\bar{\partial}\mathbb{K}_{+}(t)=\sqrt{-1}(1-\frac{s}{12}t^2)
(\varphi^2\wedge\bar{\varphi}^2\wedge\varphi^3\wedge\bar{\varphi}^3+
\varphi^3\wedge\bar{\varphi}^3\wedge\varphi^1\wedge\bar{\varphi}^1$

~~~~~~~~~~~~~~~~$-\frac{s}{24}\varphi^1\wedge\bar{\varphi}^1\wedge\varphi^2\wedge\bar{\varphi}^2),$

$d^{-1,2}d^{2,-1}\mathbb{K}_{-}(t)=\sqrt{-1}(1+\frac{s}{12}t^2)
(\varphi^2\wedge\bar{\varphi}^2\wedge\bar{\varphi}^3\wedge\varphi^3+
\bar{\varphi}^3\wedge\varphi^3\wedge\varphi^1\wedge\bar{\varphi}^1$

~~~~~~~~~~~~~~~~~~~~~~~~~$+\frac{s}{24}\varphi^1\wedge\bar{\varphi}^1\wedge\varphi^2\wedge\bar{\varphi}^2),$

\noindent where $s$ is the scalar curvature of 4-manifold $(M, g)$.

\begin{center}
\item\section{Unitary connections on complex vector bundles over the twistor space}
\end{center}
In this section, for the first Chern form of a natural unitary connection on the vertical tangent bundle over the twistor space $Z$, we can recover J. Fine and D. Panov's result on the condition of the first Chern form being symplectic and P. Gauduchon's result on the condition of the first Chern form being a (1,1)-form respectively, by using the method of moving frames.

Given an oriented Riemannian 4-manifold $(M,g)$, the Levi-Civita connection on $(M,g)$ induces
a splitting of the tangent bundle $TZ$ into the direct sum of the horizontal and vertical distributions, i.e., $TZ=\mathcal{H}\oplus \mathcal{V}$. $\mathcal{H}\cong \pi_2^\star TM$, $\mathcal{V}=\text{Ker}(d\pi_2)$, where $\pi_2$ is the projection $\pi_2:Z\rightarrow M$. In fact, from section 2, we can construct a natural Riemannian metric
$g_t$ on the twistor space, in this case, locally, $\mathcal{H}=\text{Span}_{\mathbb{R}}\{u^{\star}\theta^\alpha|\alpha=1,2,3,4\}$, and
$\mathcal{V}=\text{Span}_{\mathbb{R}}\{u^{\star}\theta^5,u^{\star}\theta^6\}$, where $u$
is any local section of projection $\pi_1:P\rightarrow Z$.

For almost Hermitian twistor spaces $(Z,g_t,\mathbb{J}_{\pm})$,
the vector bundle $\mathcal{H}$ can be considered as
 a Hermitian vector bundle of rank 2, and $\mathcal{V}$ is a Hermitian line bundle. It is well known that \cite{Hi} the first
Chern classes of complex tangent bundles $(TZ,\mathbb{J}_{\pm})$ satisfy
$$c_1(Z,\mathbb{J}_-)=0,~~c_1(Z,\mathbb{J}_+)=2c_1(\mathcal{V}).$$

We only study a natural unitary connection on complex tangent bundle $(TZ,\mathbb{J}_{+})$ with metric $g_t$, and the properties of the first Chern forms of the natural complex vector bundles $\mathcal{H}$ and $\mathcal{V}$ over the twistor space with the induced unitary connections.

With respect to complex tangent bundle $(TZ,\mathbb{J}_{+})$ with metric $g_t$, the
local unitary co-frame is $\{u^{\star}\varphi^1,u^{\star}\varphi^2,2tu^{\star}\varphi^3\}$. From the structure equations (3.1)-(3.4), we have
\begin{equation}
d\begin{pmatrix} \varphi^1\\ \varphi^2\\ 2t\varphi^3\end{pmatrix}=
\left(\begin{array}{ccc}\sqrt{-1}\omega_2^1&-\psi^3&0\\
\bar{\psi}^3&\sqrt{-1}\omega_4^3&0\\
0&0&\sqrt{-1}(\omega_2^1+\omega_4^3)\end{array}\right)\wedge \begin{pmatrix} \varphi^1\\ \varphi^2\\ 2t\varphi^3\end{pmatrix}
+\begin{pmatrix} \bar{\varphi}^2\wedge\varphi^3\\ \varphi^3\wedge\bar{\varphi}^1\\ 2t\rho(\Omega)\end{pmatrix}.
\end{equation}

Set
\begin{equation}
\omega_{\mathbb{J}_+}=-\left(\begin{array}{ccc}\sqrt{-1}\omega_2^1&-\psi^3&0\\
\bar{\psi}^3&\sqrt{-1}\omega_4^3&0\\
0&0&\sqrt{-1}(\omega_2^1+\omega_4^3)\end{array}\right),~\tau=\begin{pmatrix} \bar{\varphi}^2\wedge\varphi^3\\ \varphi^3\wedge\bar{\varphi}^1\\ 2t\rho(\Omega)\end{pmatrix}.
\end{equation}
Obviously, $\omega_{\mathbb{J}_+}+\bar{\omega}_{\mathbb{J}_+}^T=0$.

On the same open subset $U\subset Z$, if $\tilde{u}:U\rightarrow P$ is
another local section of the $U(2)$-principal bundle $\pi_1:P\rightarrow Z$, then $\tilde{u}=u\varrho(\mathfrak{a}^{-1})$, where $\varrho$ is the isomorphism in (2.11),
$\mathfrak{a}$ is a smooth $U(2)$-valued function on $U$. We can also write $\tilde{u}=R_{\varrho(\mathfrak{a}^{-1})}u$.
Here $R_{\varrho(\mathfrak{a}^{-1})}:P|_{U}\rightarrow P|_{U}$ is defined as follows:
for $\forall J_{x}\in U$, $\forall e\in P_{J_{x}}$,
$R_{\varrho(\mathfrak{a}^{-1})}(e):=e\varrho(\mathfrak{a}^{-1}(J_{x}))$. In this case,   $U(2)$-valued function $\mathfrak{a}$ may be considered as a fiber-constant function on $P|_{U}$.
From the structure equations (2.4) and (2.5), we have
\begin{eqnarray}
&&R_{\varrho(\mathfrak{a}^{-1})}^{\star}\theta=\varrho(\mathfrak{a})\theta,\\
&&R_{\varrho(\mathfrak{a}^{-1})}^{\star}\Omega=\varrho(\mathfrak{a})\Omega {\varrho(\mathfrak{a}^{-1})},\\
&&R_{\varrho(\mathfrak{a}^{-1})}^{\star}\omega=\varrho(\mathfrak{a})\omega \varrho(\mathfrak{a}^{-1})
-d\varrho(\mathfrak{a})\cdot\varrho(\mathfrak{a}^{-1}).
\end{eqnarray}

It follows that
\begin{equation}
R_{\varrho(\mathfrak{a}^{-1})}^{\star}\begin{pmatrix} \varphi^1\\ \varphi^2\\ 2t\varphi^3\end{pmatrix}=
\left(\begin{array}{cc}\mathfrak{a}&0\\
0&\text{det}(\mathfrak{a})\end{array}\right)\begin{pmatrix} \varphi^1\\ \varphi^2\\ 2t\varphi^3\end{pmatrix}:=\mathfrak{g}\begin{pmatrix} \varphi^1\\ \varphi^2\\ 2t\varphi^3\end{pmatrix},
\end{equation}
and
\begin{equation}
R_{\varrho(\mathfrak{a}^{-1})}^{\star}\rho(\Omega)=\text{det}(\mathfrak{a})\rho(\Omega),~~~
R_{\varrho(\mathfrak{a}^{-1})}^{\star}\tau=\mathfrak{g}\tau.
\end{equation}
From (4.1), (4.6) and (4.7), we obtain
\begin{equation}
\tilde{u}^{\star}\omega_{\mathbb{J}_+}=\mathfrak{g}u^{\star}\omega_{\mathbb{J}_+}\mathfrak{g}^{-1}
-d\mathfrak{g}\cdot\mathfrak{g}^{-1}.
\end{equation}
It follows that $\omega_{\mathbb{J}_+}$ can be used to define an unitary connection (also called almost Hermitian connection) on the complex tangent bundle $(TZ,\mathbb{J}_{+})$.  From (4.1), the connection $\omega_{\mathbb{J}_+}$ is not a Chern connection, because the torsion $\tau$ has components of $(1,1)$-form. Meanwhile, from $\omega_{\mathbb{J}_+}$,  we can define the induced unitary connections on Hermitian vector bundles $\mathcal{H}$ and $\mathcal{V}$, respectively. In particular, this induced unitary connection on $\mathcal{V}$ is the same as the connection in \cite{FP} or \cite{Gau}.

By the structure equations, we can calculate the curvature, denoted by $\Omega_{\mathbb{J}_+}=
d\omega_{\mathbb{J}_+}+\omega_{\mathbb{J}_+}\wedge \omega_{\mathbb{J}_+}$, with respect to connection $\omega_{\mathbb{J}_+}$. Thus the first Chern form, denoted by $c_1(Z,\omega_{\mathbb{J}_+})$, is
\begin{equation}
c_1(Z,\omega_{\mathbb{J}_+})=\frac{\sqrt{-1}}{2\pi}\text{trace}(u^{\star}\Omega_{\mathbb{J}_+})
=\frac{1}{\pi}(4u^{\star}\theta^5\wedge u^{\star}\theta^6+u^{\star}\Omega_2^1+u^{\star}\Omega_4^3).
\end{equation}

On the Hermitian vector bundles $\mathcal{H}$ and $\mathcal{V}$ with the induced unitary connections, the first Chern forms denoted by $c_1(\mathcal{H},\omega_{\mathbb{J}_+})$ and $c_1(\mathcal{V},\omega_{\mathbb{J}_+})$, respectively, are
\begin{equation}
c_1(\mathcal{H},\omega_{\mathbb{J}_+})=c_1(\mathcal{V},\omega_{\mathbb{J}_+})
=\frac{1}{2\pi}(4u^{\star}\theta^5\wedge u^{\star}\theta^6+u^{\star}\Omega_2^1+u^{\star}\Omega_4^3).
\end{equation}

Now, we can re-prove the following two interesting theorems in \cite{FP} and \cite{Gau} respectively.
\vspace{0.2cm}

\noindent\textbf{Theorem 4.1}\cite{FP}~~\emph{$c_1(\mathcal{V},\omega_{\mathbb{J}_+})$ (or $c_1(Z,\omega_{\mathbb{J}_+})$) is a symplectic form if and only if the endomorphism of the bundle of self-dual 2-forms $\wedge^+$ given by
$$\mathfrak{D}=
(W^++\frac{s}{12}\text{Id})^2-\text{Ric}_0^\star \text{Ric}_0:\wedge^+\rightarrow \wedge^+$$
is definite.}

\vspace{0.2cm}

\emph{Proof.} From (4.9) or (4.10), it follows $c_1(\mathcal{V},\omega_{\mathbb{J}_+})$ (or $c_1(Z,\omega_{\mathbb{J}_+})$) is a symplectic form if and only if
\begin{equation}
(\Omega_2^1+\Omega_4^3)\wedge(\Omega_2^1+\Omega_4^3)\neq0
\end{equation}
on principal bundle $P$.

For the given oriented Riemannian 4-manifold $(M,g)$, the Levi-Civita connection induces
a metric connection on the bundle of self-dual 2-forms $\wedge^+$, also called Levi-Civita connection. The curvature of this connection on $\wedge^+$ is denoted by $F_{\nabla}\in \Gamma(\wedge^2\otimes \mathfrak{so}(\wedge^+))$. As in \cite{AHS} or \cite{FP}, using
the identification $\wedge^+\cong \mathfrak{so}(\wedge^+)^{\star}$, we have
$$F_{\nabla}=(W^{+}+\frac{s}{12}\text{Id})\oplus Ric_{0}:~\wedge^{+}\longrightarrow \wedge^{+}\oplus \wedge^{-}.$$
Theorem 3.1 in \cite{FP} shows that the image of $F_{\nabla}$ is a maximal definite subspace of
$\wedge^2$ for the wedge product if and only if the endomorphism $\mathfrak{D}$ is definite.

Let $e=(e_{\alpha}):V\subset M\rightarrow P$ be a local section, its dual co-frame is $(e^{\star}\theta^{\alpha})$. Then $\{e^{\star}\alpha_{+}^1,e^{\star}\alpha_{+}^2,e^{\star}\alpha_{+}^3,\}$ is an orthonormal
basis for $\wedge^{+}$. With respect to this basis, by direct calculations, the curvature map $F_{\nabla}$ of the  Levi-Civita connection on $\wedge^+$ has the following matrix representation:
$$\Phi_{\wedge^{+}}=\left(\begin{array}{ccc}0&e^{\star}\Omega_3^2+e^{\star}\Omega_4^1&
e^{\star}\Omega_4^2-e^{\star}\Omega_3^1\\
-e^{\star}\Omega_3^2-e^{\star}\Omega_4^1&0&e^{\star}\Omega_2^1+e^{\star}\Omega_4^3\\
e^{\star}\Omega_3^1-e^{\star}\Omega_4^2&-e^{\star}\Omega_2^1-e^{\star}\Omega_4^3&0\end{array}\right).$$

Since the local section $e=(e_{\alpha}):V\subset M\rightarrow P$ is arbitrary, it is obvious that
the image of $F_{\nabla}$ is a maximal definite subspace of
$\wedge^2$ if and only if $(\Omega_2^1+\Omega_4^3)\wedge(\Omega_2^1+\Omega_4^3)\neq0$
on principal bundle $P$.
\hfill$\Box$

\vspace{0.2cm}

\noindent\textbf{Theorem 4.2}\cite{Gau}~~\emph{$c_1(\mathcal{V},\omega_{\mathbb{J}_+})$ (or $c_1(Z,\omega_{\mathbb{J}_+})$) is a $(1,1)$-form if and only if the given oriented Riemannian 4-manifold $(M,g)$ is an anti-self-dual 4-manifold.}

\vspace{0.2cm}

\emph{Proof.} From the definition of almost complex structure $\mathbb{J}_+$ and (4.9)-(4.10),
$c_1(\mathcal{V},\omega_{\mathbb{J}_+})$ (or $c_1(Z,\omega_{\mathbb{J}_+})$) is a $(1,1)$-form if and only if
\begin{equation}
\Omega_2^1+\Omega_4^3\equiv0~\text{mod}~(\varphi^{i}\wedge\bar{\varphi}^{j},~~i,j=1,~2)
\end{equation}
on principal bundle $P$.

Since
\begin{equation}
\Omega_2^1+\Omega_4^3=\frac{\sqrt{-1}}{2}[(R_{1212}+R_{3412})\varphi^1\wedge\bar{\varphi}^1
+(R_{3434}+R_{1234})\varphi^2\wedge\bar{\varphi}^2]
\end{equation}

~~~~~~~~~~~~~~$+\frac{1}{2t^2}(\Lambda_2\varphi^1\wedge\bar{\varphi}^2
+\bar{\Lambda}_2\bar{\varphi}^1\wedge\varphi^2+\Lambda_3\varphi^1\wedge\varphi^2
+\bar{\Lambda}_3\bar{\varphi}^1\wedge\bar{\varphi}^2)$,

\noindent where
$$\Lambda_3=\frac{t^2}{2}[R_{1213}+R_{3413}-R_{1224}-R_{3424}
-\sqrt{-1}(R_{1214}+R_{3414}+R_{1223}+R_{3423})].$$

From (4.13), condition (4.12) is equivalent to $\Lambda_3=0$ on principal bundle $P$. As
$$A_{12}=\frac{1}{2}(R_{1213}+R_{3413}-R_{1224}-R_{3424}),
~A_{13}=\frac{1}{2}(R_{1214}+R_{3414}+R_{1223}+R_{3423}),$$
then condition (4.12) is equivalent to $A_{12}=A_{13}=0$ on principal bundle $P$.

The method is the same as the proof of Theorem 3.1, it follows that $A_{12}=A_{13}=0$ on principal bundle $P$ if and only if the given oriented Riemannian 4-manifold $(M,g)$ is an anti-self-dual 4-manifold.
\hfill$\Box$

\vspace{0.2cm}

\noindent\textbf{Remark 4.3}~~As an example, for an anti-self-dual Einstein 4-manifold $(M, g)$ with scalar curvature $s\neq0$, from (4.13), we have
$$\Omega_2^1+\Omega_4^3=\frac{\sqrt{-1}}{24}s(\varphi^1\wedge\bar{\varphi}^1
+\varphi^2\wedge\bar{\varphi}^2),$$
$$(\Omega_2^1+\Omega_4^3)\wedge(\Omega_2^1+\Omega_4^3)=\frac{s^2}{72}
\theta^1\wedge\theta^2\wedge\theta^3\wedge\theta^4,$$
on principal bundle $P$.
We can also define an unitary connection $\omega_{\mathbb{J}_-}$ on the complex tangent bundle $(TZ,\mathbb{J}_{-})$ with metric $g_t$ as in (4.1), and we also obtain   similar results as in Theorem 4.1 and Theorem 4.2. In fact, by direct calculations, the first Chern forms are given by
$$c_1(Z,\omega_{\mathbb{J}_-})=0,~~
c_1(\mathcal{H},\omega_{\mathbb{J}_-})=-c_1(\mathcal{V},\omega_{\mathbb{J}_-})
=\frac{1}{2\pi}(4u^{\star}\theta^5\wedge u^{\star}\theta^6+u^{\star}\Omega_2^1+u^{\star}\Omega_4^3).$$

\vspace{0.3cm}
\noindent\textbf{Acknowledgments}\\
\noindent Zhou would like to thank Professor Jiagui Peng  and Professor Xiaoxiang Jiao for their helpful suggestions and encouragements. Fu is supported in part by NSFC grants
10831008 and 11025103.

\end{document}